\documentclass[reqno,11pt]{amsart}

\usepackage[a4paper,  margin=3cm]{geometry}

\usepackage{amsmath,amssymb}
\usepackage{dsfont}         
\usepackage{comment}
\newtheorem{theorem}{Theorem}[section]
\newtheorem{corollary}[theorem]{Corollary}
\newtheorem{lemma}[theorem]{Lemma}
\newtheorem{proposition}[theorem]{Proposition}
\theoremstyle{definition}
\newtheorem{definition}[theorem]{Definition}
\newtheorem{assumption}[theorem]{Assumption}
\theoremstyle{remark}
\newtheorem{remark}[theorem]{Remark}
\numberwithin{equation}{section}

\title[States of an infinite system with branching]{Evolution of states of an infinite particle system with nonlocal branching}

\author{Yuri Kozitsky}

\address{Instytut Matematyki, Uniwersytet Marii Curie-Sk{\l}odowskiej, 20-031 Lublin, Poland}
\email{jkozi@hektor.umcs.lublin.pl}

\author{Agnieszka Tana\'s}

\address{Politechnika Lubelska, 20-618 Lublin, Poland}


\begin{document}
\maketitle

\begin{abstract}
We study the evolution of states of an infinite system of point
particles dwelling in a locally compact Polish space $X$. Each
particle produces at random  a finite `cloud' of offsprings
distributed over $X$ according to some law, and disappears
afterwards. The system's states are probability measures on an
appropriate space of locally finite counting measures on $X$.  Their
evolution is obtained by solving the corresponding Fokker-Planck
equation. We prove that this equation has a unique solution and
discuss some of its properties. Our pivotal idea of dealing with
infinite systems consists in passing to tempered counting measures
by imposing appropriate restrictions on the branching. In this
approach, we first solve a nonlinear evolution equation in the space
of bounded continuous functions on $X$ -- so called log-Laplace
equation. Next we solve the Kolmogorov equation which is then used
to solve the Fokker-Planck equation and thus describe the evolution
in question.
\end{abstract}

\section{Introduction}

In recent years, the stochastic evolution of infinite particle
systems attract considerable attention, see, e.g.,
\cite{Konar,Konar1,KP,Koz2}. A related popular topic is
measure-valued stochastic branching characterizing the evolution of
random `clouds' \cite{Dawson,DKS,Li}, see also
\cite{BL,BL1,DGL,Koz1} and the literature quoted in these works. Let
$X$ be a locally compact Polish space, $\mathcal{B}(X)$ its Borel
$\sigma$-field and $\mathcal{N}$ be the set of all finite counting
measures on $X$, i.e., $\nu (\Delta)$ is a nonnegative integer for
each $\nu \in \mathcal{N}$ and  $\Delta \in \mathcal{B}(X)$. It is
known that the weak topology makes $\mathcal{N}$ a Polish space. By
\cite[Proposition 9.1.III, page 4]{DV2} it follows that $\nu$ can be
presented in the form $\nu = \sum_i  \delta_{x_i}$, where
$\delta_{x_i}$ are Dirac's measures and some of $x$'s may coincide.
By this formula one may interpret $\nu$ as a `cloud' of particles
located at points $x_i\in X$. The key aspect of this interpretation
is that $\nu(X)$ -- the total number of particles -- is finite as
$\nu$ is a finite measure. Since in the course of branching each
particle is replaced by a finite number of offsprings, the system
remains finite during all its lifetime. Our main challenge in this
work is to find out how to describe branching in an infinite
particle system, in which the number of offsprings appearing at a
given moment of time can be infinite. Along with a purely
mathematical meaning of this question, its application value is
related to the known fact that the collective behavior of a
macroscopically large system can be understood only in the `infinite
system limit', see, e.g., \cite[pages 5,6]{Simon}.

The way of dealing with an infinite particle system which we propose
in this work has the following main aspects. Our object is an
infinite collection of branching point particles -- an infinite
cloud -- placed in a locally compact Polish space $X$ in such a way
that each compact $\Lambda\subset X$ contains only finitely many
elements of the cloud. This means that the corresponding counting
measure belongs to the set of all locally finite counting measures
$\mathcal{N}^\#$ and may take infinite values. The branching
mechanism is described by a probability kernel $b$, i.e., a map $(X,
\mathcal{B}(\mathcal{N}))\ni (x,\Xi) \mapsto b_x(\Xi)\in [0,1]$ such
that each $b_x$ is a probability measure on $\mathcal{N}$ and
$x\mapsto b_x(\Xi)$ is measurable for each  Borel subset $\Xi$ of
the set of all finite counting measures $\mathcal{N}$. Let
$\delta(x)$ be the probability that a point at $x$  disappears
without leaving offsprings. That is, $\delta(x) = b_x(\Gamma^0)$,
where $\Gamma^0$ is the singleton consisting of the zero measure.
Our pivotal idea is to impose the condition that $1-\delta(x)$
vanishes at infinity, i.e., $1-\delta(x)< \varepsilon$ whenever
$x\in \Lambda^c_\varepsilon:=X\setminus \Lambda_\varepsilon$, for a
sufficiently big compact $\Lambda_\varepsilon\subset X$. Then we
consider only those $\nu\in \mathcal{N}^\#$ -- `tempered measures'
-- for which $1-\delta(x)$ is integrable. Note that imposing a
condition of this kind seems inevitable as an infinite system of
branching particles can produce simultaneously an infinite cloud of
offsprings that collapses into a compact $\Lambda\subset X$, and
thus destroys the aforementioned local finiteness of the cloud. The
same problem arises also in the dynamical theory of infinite systems
of physical particles, see \cite[page 223]{DSS}, where it is settled
by imposing similar restrictions.

In dealing with particle systems, it is more convenient for us to
stick at the `corpuscular' terminology, i.e., to speak of {\em
configurations} of particles instead of counting measures. Following
\cite{Lenard}, by a configuration $\gamma$ we mean a countable
collection of point particles placed in $X$, where each particle is
completely characterized by its location $x\in X$. Multiple
locations are possible and each compact $\Lambda\subset X$ may
contain only finitely many elements of $\gamma$. That is,
`configuration' is a more rigorous synonym of the aforementioned
`cloud'. The set of all configurations is denoted by $\Gamma$. Note
that particles with the same location are indistinguishable, and
there can only be finitely many of them located at a given $x$. By
writing $\gamma\cup x$ we mean the configuration with added particle
located at $x$. Likewise we define $\gamma\setminus x$ for $x\in
\gamma$. Then by $\sum_{x\in \gamma}$ we mean $\sum_{i}$ for a
certain enumeration of the elements of $\gamma$, cf. \cite{Lenard}.
In this context, each $\nu \in \mathcal{N}^\#$ is presented as
$\sum_{x\in \gamma} \delta_x$, which establishes a bijection between
$\Gamma$ and $\mathcal{N}^\#$, see above. Typically,
$\mathcal{N}^\#$ is equipped with the vague (weak-hash) topology
which is the weakest topology that makes continuous the maps $\nu
\mapsto \int_X g d \nu =:\nu(g)$ with all choices of compactly
supported continuous functions $g:X\to \mathds{R}$. Then the same
topology on $\Gamma$ is defined by the maps $\gamma \mapsto
\sum_{x\in \gamma} g(x)$. This makes $\Gamma$ and $\mathcal{N}^\#$
Polish spaces, see \cite[Proposition 9.1,IV, page 6]{DV2}.

Our model is defined by the Kolmogorov operator
\begin{equation}
  \label{6F}
  (LF) (\gamma) =  \sum_{x\in \gamma} \int_{\Gamma}\left[F(\gamma\setminus x\cup \xi) - F(\gamma)
  \right] b_x( d \xi),
\end{equation}
where $F$ is a suitable (test) function and $b$ is the
aforementioned branching kernel. Note that the sum in (\ref{6F}) is
infinite for infinite $\gamma$. The Kolmogorov operator (\ref{6F})
describes the distribution of the offsprings (constituting cloud
$\xi$) of a particle located at $x$. Its detailed properties are
listed in Assumption \ref{Ass1} below. By means of this $L$ we
introduce the Kolmogorov equation
\begin{equation}
  \label{KE}
  \frac{d}{dt} F_t = L F_t, \qquad F_{t}|_{t=0} = F_0,
\end{equation}
describing the evolution of test functions, and the Fokker-Planck
equation
\begin{equation}
  \label{FPE}
  \mu_t (F) = \mu_0 (F) + \int_0^t \mu_s ( LF) d s,
\end{equation}
that describes the evolution of states of the considered system.
Here the states are probability measures on $\Gamma$; the set of all
such states is denoted by $\mathcal{P}(\Gamma)$. Finally,  $\mu(F):=
\int F d \mu$ for suitable $F:\Gamma\to \mathds{R}$.  A
comprehensive theory of the evolution equations of this kind can be
found in \cite{Mich}.

The main steps in solving (\ref{FPE}) undertaken in this work can be
outlined as follows. In view of the form of $L$ given in (\ref{6F}),
a natural way of solving (\ref{KE}) is to assume that $F_t$ has the
following form
\begin{equation}
  \label{Gu}
  F_t (\gamma) = \prod_{x\in \gamma} \phi_t(x),
\end{equation}
where, for each fixed $t$, $\phi_t$ is a continuous functions of
$x\in X$ taking values in $(0,1)$. It turns out that the Kolmogorov
equation admits solutions of this type with $\phi_t$ satisfying a
nonlinear evolution equation derived from (\ref{KE}). In the theory
of branching processes, such equations are called `log-Laplace', see
\cite[pages 60, 61]{Dawson}. We show that the latter has a unique
solution and thereby obtain solutions of (\ref{KE}) in the Banach
space of bounded continuous functions $F$, where $L$ is defined as a
closed linear operator with a suitable domain $\mathcal{D}(L)$.
Having this done, we show that (\ref{FPE}) has a solution $t \mapsto
\mu_t$ with $F\in \mathcal{D}(L)$. It turns out that the latter is
big enough so that the following holds: (a) for each bounded
measurable function $F'$, $\mu(F')$ can be approximated by $\mu(F)$,
$F\in \mathcal{D}(L)$; (b) the mentioned solution be unique. The
main result of this work is Theorem \ref{2tm} that states that
(\ref{FPE}) has a unique solution possessing a number of properties
presented in this statement. This result is based on Lemma \ref{1tm}
that states the existence of classical solutions of the Kolmogorov
equation (\ref{KE}) in the form of (\ref{Gu}), with $t \mapsto
\phi_t$ described in Lemma \ref{LUlm}.

The rest of this paper has the following structure. In Sect.2, we
settle the mathematical framework and formulate our assumptions
concerning the branching kernel $b$. Then we introduce and describe
tempered configurations by employing a continuous function $\psi:X
\to (0,+\infty)$, that vanishes at infinite and is such that the
aforementioned death probability satisfies $\delta(x)\geq 1 -
\psi(x)$. The key statement of this part is Proposition \ref{F1pn}
according to which the set of all tempered configurations
$\Gamma^\psi$ is a Polish space. This allows us to restrict
ourselves to considering the states with the property
$\mu(\Gamma^\psi)=1$. In fact, this restriction is a direct analog
of the condition imposed in \cite{DSS}, see {\it ibid,} Definitions
3.1 and 3.2, page 223. Next, we discuss in detail the properties of
the branching kernel. In Sect. 3, we prepare solving our main
evolution equation (\ref{FPE}) by defining $L$ as a closed linear
operator with domain $\mathcal{D}(L)$ in the Banach space of bounded
continuous functions $F:\Gamma^\psi \to \mathds{R}$. As mentioned
above, the key ingredient of this construction is solving the
log-Laplace lequation in the space of bounded continuous functions
$\phi:X \to \mathds{R}$, defined by the branching kernel, see Lemma
\ref{LUlm}. This step is typical in the theory of branching
processes, cf. \cite[Theorem 3.1]{BL1}. With the help of this, we
prove that the Kolmogorov equation (\ref{KE}) with the just
mentioned closed Kolmogorov operator  $L$  has a classical solution,
see Lemma \ref{1tm}, in the form of (\ref{Gu}) with $\phi_t$ solving
the log-Laplace equation. Thereby and after additional preparations,
we prove (Theorem \ref{2tm}) that (\ref{FPE}) has a unique solution
$t\mapsto \mu_t$, which is weakly continuous, i.e., $\mu_t
\Rightarrow \mu_s$ as $t\to s$. In the subcritical case, we show
that $\mu_t \Rightarrow \mu_\infty$ as $t\to +\infty$, where
$\mu_\infty(\Gamma^0)=1$. At the very end, we make some concluding
remarks on possible extensions of the results of this work.

\section{Preliminaries and Assumptions}

\subsection{Notions and notations}
By $\mathds{1}_A$ we denote the indicator of a suitable set $A$. A
Polish space is a separable topological space that can be metrized
by a complete metric, see \cite[Chapt. 8]{Cohn}. For a Polish space
$E$, $\mathcal{B}(E)$ will stand for the corresponding Borel
$\sigma$-field. By $C_{\rm b}(E)$, $C_{\rm cs}(E)$, $B_{\rm b}(E)$
we denote the sets of all continuous and bounded, continuous and
compactly supported, measurable and bounded functions $f:E\to
\mathds{R}$, respectively. By $C_{\rm b}^{+}(E)$ we denote the set
of positive elements of $C_{\rm b}(E)$. Finally, by $C_{0}^{+}(E)$
we denote the set of all $f\in C_{\rm b}^{+}(E)$ which satisfy: (a)
$f(x) >0$ for all $X$; (b) for each $\varepsilon
>0$, one finds a compact $\Lambda_\varepsilon \subset X$ such that
$f(x) < \varepsilon$ whenever $x\in X\setminus \Lambda_\varepsilon$.

Let $\mathcal{F}$ be a family of functions $f:E\to \mathds{R}$. By
$\sigma \mathcal{F}$ we denote the smallest sub-field of
$\mathcal{B}(E)$ such that each $f\in \mathcal{F}$ is $\sigma
\mathcal{F}/\mathcal{B}(\mathds{R})$-measurable. By $\mathcal{P}(E)$
we denote the set of all probability measures on
$(E,\mathcal{B}(E))$; for suitable $f:E\to \mathds{R}$ and $\mu\in
\mathcal{P}(E)$, we write $\mu(f) = \int_E f d \mu$. The weak
topology of $\mathcal{P}(E)$ is defined as the weakest one that
makes continuous all the maps $\mu \mapsto \mu(f)$, $f\in C_{\rm
b}(E)$. With this topology $\mathcal{P}(E)$ is also a Polish space.
By writing $\mu_n \Rightarrow \mu$, $n\to +\infty$, we mean that
$\{\mu_n\}_{n\in \mathds{N}}$ weakly converges to $\mu$. A family
$\mathcal{F}$ of functions $f:E\to \mathds{R}$ is called
\emph{separating} if $\mu_1(f) = \mu_2(f)$, holding for all $f\in
\mathcal{F}$, implies $\mu_1 = \mu_2$ for each pair $\mu_1,\mu_2\in
\mathcal{P}(E)$. Furthermore, $\mathcal{F}$ is said to separate the
points of $E$ if for each distinct $x,y\in E$, one finds
$f\in\mathcal{F}$ with the property $f(x)\neq f(y)$. If
$\mathcal{F}$ separates points and is closed with respect to
multiplication, it is separating, see \cite[Theorem 4.5, page
113]{EK}. A family $\mathcal{F}$ is called \emph{convergence
determining} if $\mu_n (f) \to \mu(f)$, holding for all $f\in
\mathcal{F}$, implies $\mu_n \Rightarrow \mu$.
\subsection{Tempered configurations}

As mentioned above, we will deal with a locally compact Polish space
$X$. Usually, branching processes are constructed in a more general
setting, cf. \cite{Li}, e.g., by taking as $X$ a Luzin space
\cite{BL,BL1}. Our choice of $X$ was done mostly for the following
reason. In contrast to \cite{BL,BL1}, we study branching in infinite
systems and thus deal with the space of infinite configurations. In
Polish spaces without assuming their local compactness, the vague
topology of spaces of infinite configurations is introduced with the
help of continuous functions vanishing outside bounded sets, cf.
\cite[page 2]{DV2}, which assumes fixing some concrete metric of
$X$. Thus, in order to be free in choosing such metrics, as well as
to avoid possible complications of the topological aspects of our
work, we restrict ourselves here to considering locally compact
spaces. Note that the local compactness together with the
separability of $X$ imply its $\sigma$-compactness, i.e, the
existence of a nest of compact subsets $\{\Lambda_k\}_{k\in
\mathds{N}}$ that exhaust $X$. The existence of such nests and the
corresponding compact Lyapunov functions could be taken as the basic
topological assumption concerning $X$, as it was done in
constructing quasi-regular Dirichlet forms \cite{MaR}, or cadlag
processes by means of resolvent kernels \cite{BeR}.

Among all infinite configurations, one may distinguish those that
have a priori prescribed properties. Here we do this by employing a
function $\psi\in C^{+}_{\rm b}(X)$, $\psi(x) \leq 1$, for which we
set
\begin{equation}
\label{7F}
 \Psi (\gamma) =
\sum_{x\in \gamma}\psi(x).
\end{equation}
Then we define the set of tempered configurations as
\begin{equation*}
  \Gamma^\psi =\{\gamma\in \Gamma: \Psi(\gamma) < \infty\}.
\end{equation*}
It is clear that
\begin{equation}
  \label{8Fa}
\Gamma^{\psi'} \supset \Gamma^\psi, \qquad {\rm whenever} \ \
\psi'\leq \psi.
\end{equation}
By this observation we can vary $\Gamma^\psi$ from $\Gamma$ (by
taking $\psi\in C_{\rm cs}^{+}(X)$) to $\Gamma_0:=\{\gamma\in
\Gamma: \gamma \ {\rm is} \ {\rm finite}\}$, corresponding to
$\psi\equiv 1$. If $\psi\in C_0^{+}(X)$, then $\Gamma^\psi$ is a
proper subset of $\Gamma$ and supset of $\Gamma_0$. As an example.
one can take $X=\mathds{R}$ and $\psi(x) = \psi_0 e^{-\alpha |x|}$,
$\alpha
>0$, $\psi_0\in (0,1]$. Then the configuration $\mathds{N}\subset \mathds{R}$ is in
$\Gamma^{\psi}$, whereas $\{ \log n: n\in \mathds{N}\}$ is not if
$\alpha \leq 1$. In a separate publication, we plan to study scales
of such spaces $\Gamma^\psi$, including the inductive and projective
limit topologies on $\cup_{\psi\in \varPsi}\Gamma^\psi$ and
$\cap_{\psi\in \varPsi'}\Gamma^\psi$, respectively, with suitable
families $\varPsi$, $\varPsi'$.

In the sequel, we employ one and the same $\psi\in C_0^{+}(X)$,
separated away from zero, i.e., such that $\inf_{x\in \Lambda}
\psi(x) > 0$ for each compact $\Lambda \subset X$. Its  choice will
be done in the next subsection.  For each $\gamma\in \Gamma^\psi$,
the measure
\begin{equation}
  \label{9F}
  \nu_\gamma = \sum_{x\in\gamma} \psi(x) \delta_x
\end{equation}
is finite. Thus, one can equip $\Gamma^\psi$ with the topology
defined as the weakest one that makes continuous all the maps
\begin{equation}
  \label{10F}
\Gamma^\psi \ni \gamma \mapsto \sum_{x\in \gamma} g(x) \psi(x),
\qquad g\in C_{\rm b}(X).
\end{equation}
Similarly as in Proposition 2.7 and Corollary 2.8 of \cite{Koz2}, we
prove the following.
\begin{proposition}
  \label{F1pn}
With the topology defined in (\ref{10F}), $\Gamma^\psi$  is a Polish
space, continuously embedded in $\Gamma$. Thus,
$\mathcal{B}(\Gamma^\psi) =\{ A\in \mathcal{B}(\Gamma): A \subset
\Gamma^\psi\}$.
\end{proposition}
\begin{proof}
First we note that the set of measures $\{\nu_\gamma:\gamma\in
\Gamma^\psi\}$ is a subset of the space $\mathcal{M}$ of all finite
positive Borel measures on $X$, which is a Polish space with the
weak topology. Let us prove that
$\{\nu_\gamma:\gamma\in\Gamma^\psi\}$ is a closed subset of
$\mathcal{M}$. To this end, we take a sequence $\{\gamma_n\}_{n\in
\mathds{N}}\subset \Gamma^\psi$ such that $\{\nu_{\gamma_n}\}_{n\in
\mathds{N}}$ is a Cauchy sequence in a metric of $\mathcal{M}$ that
makes this space complete. Let $\nu \in \mathcal{M}$ be its limit,
and hence
\begin{equation}
  \label{Lu}
  \sum_{x\in \gamma_n} g(x) \psi(x) \to \nu (g), \qquad n \to
  +\infty,
\end{equation}
holding for all $g\in C_{\rm b}(X)$, in particular for $g\in C_{\rm
cs}(X)$ . Since $\psi$ is separated away from zero, each $h\in
C_{\rm cs}(X)$ can be written in the form $h(x) = g(x) \psi (x)$
with $g\in C_{\rm cs}(X)$. It is known, see, e.g., \cite[page
397]{Zessin}, that there exists a countable family $\{h_k\}_{k\in
\mathds{N}}\subset C_{\rm cs}(X)$ such that
\[
\bar{\upsilon}(\gamma, \gamma') := \sum_{k\in \mathds{N}} 2^{-k}
\frac{\left| \sum_{x\in \gamma}  h_k(x) - \sum_{x\in \gamma'}
h_k(x)\right|}{1+ \left| \sum_{x\in \gamma}  h_k(x) - \sum_{x\in
\gamma'} h_k(x)\right|}
\]
is a complete metric of $\Gamma$. Then the convergence as in
(\ref{Lu}) yields that the sequence $\{\gamma_n\}_{n\in \mathds{N}}$
converges to some $\gamma \in \Gamma$ in the vague topology of
$\Gamma$. To prove that this $\gamma$ lies in $\Gamma^\psi$, we take
an ascending sequence of compact $\Lambda_m \subset X$, i.e., such
that each $\Lambda_m$ lies in the interior of $\Lambda_{m+1}$ and
each $x\in X$ is contained in some $\Lambda_m$. Then we take
$g^{(m)}\in C_{\rm cs}(X)$ such that $g^{(m)}(x) =1$ for $x\in
\Lambda_m$, and $g^{(m)}(x) =0$ for $x\in X\setminus \Lambda_{m+1}$,
which is possible by Urysohn's lemma. Then
\[
\sum_{x\in \gamma} g^{(m)}(x)\psi(x) = \nu(g^{(m)}) \leq \nu(X).
\]
Now we pass here to the limit $m\to +\infty$ and obtain (by the
Beppo Levi theorem) that $\Psi(\gamma) \leq \nu(X)$, which yields,
$\gamma\in \Gamma^\psi$. Thus, $\{\nu_\gamma: \gamma\in
\Gamma^\psi\}$ is closed in $\mathcal{M}$, and thereby is also
Polish, see \cite[Proposition 8.1.2, page 240]{Cohn}. This yields
the first half of the statement. The stated continuity of the
embedding $\Gamma^\psi \hookrightarrow \Gamma$ is immediate. Then
the conclusion concerning the $\sigma$-fields follows by
Kuratowski's theorem, see \cite[Theorem 3.9, page 21]{Part}.
\end{proof}
\begin{remark}
  \label{Lurk}
The continuity of the embedding $\Gamma^\psi \hookrightarrow \Gamma$
allows one to establish the following fact:
\begin{equation}
  \label{11F}
  \mathcal{P}(\Gamma^\psi) = \{ \mu \in \mathcal{P}(\Gamma):
  \mu(\Gamma^\psi)=1\}.
\end{equation}
That is, each $\mu \in \mathcal{P}(\Gamma)$ possessing the property
$\mu(\Gamma^\psi)=1$ can be redefined as a probability measure on
$\Gamma^\psi$. Therefore, by restricting ourselves to tempered
configurations -- members of $\Gamma^\psi$ -- we exclude from our
consideration all those $\mu \in \mathcal{P}(\Gamma)$ that fail to
satisfy the mentioned support condition. At the same time, the map
$\Gamma^\psi \ni \gamma \mapsto \nu_\gamma \in \mathcal{M}$, see
(\ref{9F}), defines a natural embedding $\Gamma^\psi \hookrightarrow
\mathcal{M}$ and thus pushes each $\mu\in \mathcal{P}(\Gamma^\psi)$
to a probability measure on $\mathcal{M}$. For obvious reasons,
below we do not distinguish between the elements of
$\mathcal{P}(\Gamma^\psi)$ and the corresponding push-forward
measures on $\mathcal{M}$.
\end{remark}
Let $E$ be a Polish space. Following \cite[page 111]{EK}, we say
that a sequence $\{h_n\}_{n\in \mathds{N}}\subset B_{\rm b}(E)$
converges to a certain $h\in B_{\rm b}(E)$ \emph{boundedly} and
\emph{pointwise} if: (a) $\sup_{n}\|h_n\| < \infty$; (b) $h_n(x) \to
h(x)$ for each $x\in E$. In this case, we write $h_n
\stackrel{bp}{\to} h$. A subset, $\mathcal{H}\subset B_{\rm b}(E)$,
is said to be $bp$-closed, if $\{h_n\}\subset \mathcal{H}$ and $h_n
\stackrel{bp}{\to} h$ imply $h\in \mathcal{H}$. The $bp$-closure of
$\mathcal{H}\subset B_{\rm b}(E)$ is the smallest $bp$-closed subset
of $B_{\rm b}(E)$ that contains $\mathcal{H}$. An $\mathcal{H}'$ is
$bp$-dense in $\mathcal{H}$, if the latter is the smallest
$bp$-closed set that contains $\mathcal{H}'$. The following is
known, see \cite[Proposition 4.2, page 111]{EK} and/or \cite[Lemmas
3.2.1, 3.2.3, pages 41, 42]{Dawson}.
\begin{proposition}
  \label{N1pn}
For each Polish space $E$, there exists a countable family
$\mathcal{H} \subset C_{\rm b}^{+}(E)$ that has the following
properties:  (a) the linear span of $\mathcal{H}$ is $bp$-dense in
$B_{\rm b}(E)$; (b) $\mathcal{B}(E) = \sigma\mathcal{H}$; (c) it
contains the unit function $u(x)\equiv 1$ and is closed with respect
to addition; (d)  it is separating; (e) it is convergence
determining.
\end{proposition}
Let $\mathcal{V}=\{v_l\}_{l\in \mathds{N}}\subset C_{\rm b}^{+} (X)$
be a family of functions with the property as in Proposition
\ref{N1pn}. We may and will assume that each $v_l\in \mathcal{V}$
satisfies $\inf_{X}v_l (\hat{x}) \geq c_{0,l}
>0$ for an appropriate $c_{0,l}$, cf. \cite[Remark 3.2.3, page
42]{Dawson}. Indeed, if this is not the case, instead of $v_l$ one
can take $\tilde{v}_l:=v_l + c_{0,l}$. Then the family
$\{\tilde{v}_l\}_{l\in \mathds{N}}$ has all the properties we need.
For $\gamma\in \Gamma^\psi$, we have, cf. (\ref{9F}), $\nu_\gamma
(v_l) = \sum_{x\in \gamma} v_l (x) \psi(x)$. The topology mentioned
in Proposition \ref{F1pn} is metrizable with the metric
\begin{equation}
  \label{N3}
  \upsilon_* (\gamma,\gamma') =
\sum_{l=0}^\infty \frac{2^{-l} \left|\nu_{\gamma}(v_l)-
\nu_{\gamma'}(v_l) \right|}{1+ \left|\nu_{\gamma}(v_l)-
\nu_{\gamma'}(v_l)\right|}.
\end{equation}
For $\mu \in \mathcal{P}(\Gamma^\psi)$, its Laplace transform is
defined by the expression
\begin{gather}
  \label{N6}
  \mathfrak{L}_\mu (g)  = \mu (G^g),  \qquad g \in  C^{+}_{\rm b}(X) \\[.2cm] \nonumber
G^g (\gamma) :=  \exp \left( - \nu_{\gamma}
  (g)\right) = \exp \left( - \sum_{x\in
  \gamma} g (x) \psi(x)
  \right).
\end{gather}
It is known, see \cite[Lemma 3.2.5 and Theorem 3.2.6, page
43]{Dawson}, that the Laplace transforms of probability measures on
$\mathcal{M}$ have a number of useful properties which we are going
to exploit. In view of the embedding mentioned in the second part of
Remark \ref{Lurk}, we can attribute these properties also to the
transform defined in (\ref{N6}).
\begin{proposition}
  \label{N2pn}
Let $\mathcal{V}$ be the family of functions used in (\ref{N3}).
Then:
\begin{itemize}
  \item[(i)] $\mathcal{B}(\Gamma^\psi) = \sigma \{G^v: v\in
  \mathcal{V}\}$;
  \item[(ii)] $B_{\rm b}(\Gamma^\psi)$ is the $bp$-closure of
  the linear span of
  $\{G^v: v\in
  \mathcal{V}\}$;
\item[(iii)] $\{G^v: v\in
  \mathcal{V}\}$ is separating; \item[(iv)] $\{G^v: v\in
  \mathcal{V}\}$ is convergence determining.
\end{itemize}
\end{proposition}
The proof of claim (iv) is essentially based on the concrete choice
of the metric (\ref{N3}), by which one shows that the family $\{G^v:
v\in \mathcal{V}\}$ is strongly separating, cf. \cite[page 113]{EK}.

In the sequel, we will use the functions
\begin{equation}
  \label{13F}
  \phi (x) = 1-\theta(x) =\exp\left( - g (x) \psi (x)\right),
\end{equation}
with $g\in C_{\rm b}^{+}(X)$, that includes also the choice $g\in
\mathcal{V}$.

\subsection{The branching kernel}
We assume that, for each $x\in X$, $b_x\in \mathcal{P}(\Gamma)$ is
such that $b_x(\Gamma_0)=1$, i.e., each cloud of offsprings is
finite. Recall that $\Gamma_0\in \mathcal{B}(\Gamma)$ is the set of
all finite configurations.  For $n\in \mathds{N}_0$, we set
$\Gamma^n= \{\xi\in \Gamma_0: |\xi|=n\}$. Then $b_x (\Gamma^n)$ is
the probability of producing $n$ offsprings by the particle located
at $x$. Note that $\delta(x):=b_x(\Gamma^0)$ is just the death
probability, and
\begin{equation}
  \label{19F}
  n(x) := \int_{\Gamma_0} |\xi| b_x(d \xi) =\sum_{n=1}^{\infty} n
  b_x(\Gamma_n) = \beta_x^{(1)}(X)
\end{equation}
is the expected number of offsprings of the particle located at $x$.
Here $\beta_x^{(1)}$ is the first correlation measure of $b_x$, see,
e.g., \cite{Lenard} for a rigorous definition. Here we use its
property
\begin{equation}
  \label{19Fa}
  \int_{\Gamma} \left(\sum_{y\in \xi} h(y) \right) b_x(d\xi) = \int_X h(y)
  \beta_x^{(1)} ( d y).
\end{equation}
For $\phi$ as in (\ref{13F}), we define
\begin{equation}
  \label{21F}
  (\Phi \phi)(x) = \int_{\Gamma} \left( \prod_{y\in
  \xi}\phi(y)\right)b_x ( d \xi),
\end{equation}
with the convention that $\prod_{x\in \varnothing }\phi(x) =1$.
Clearly, $0\leq (\Phi\phi)(x) \leq 1$ for each $x\in X$. Recall that
we use  $\psi$ in (\ref{7F}) in defining tempered configurations.
\begin{assumption} The probability kernel $b$ is subject to the
following conditions:
  \label{Ass1}
  \begin{itemize}
    \item[(i)] $\Phi\phi  \in C_{\rm b}(X)$ for each $\phi$ as in (\ref{13F});
    \item[(ii)] $\sup_{x\in X} n(x) =:n_* <\infty$;
\item[(iii)] the death probability $\delta$ satisfies $\delta (x) \geq 1 - \psi(x)\geq \delta_* >0$,
holding for all $x\in X$;
\item[(iv)] there exists $m>0$ such that, for all $x\in X$, the
following holds
\begin{equation}
  \label{Q}
  \int_X \psi(y) \beta^{(1)}_x ( d y) \leq n(x) m \psi(x).
\end{equation}
  \end{itemize}
\end{assumption}
By (\ref{13F}), (\ref{19Fa}), (\ref{Q}) and Jensen's inequality one
gets
\begin{gather*}
 - \log (\Phi \phi)(x) \leq \int_{\Gamma_0}
\left(-\log \prod_{y\in \xi} \phi(y) \right) b_x ( d \xi) \\[.2cm] \nonumber
 =\int_X g(y)  \psi(y) \beta^{(1)}_x( d y)  \leq \left(\sup_{x\in X}g(x)\right)  n(x) m
\psi(x).
\end{gather*}
Note that by (\ref{N6}), (\ref{21F}) and our choice $\phi =
e^{-g\psi}$ it follows that
\begin{equation*}
  (\Phi \phi)(x) = \int_{\Gamma_0} G^g (\xi) b_x ( d \xi) =
  \mathfrak{L}_{b_x}(g).
\end{equation*}
Then assumption (i) can be reformulated as the continuity of the map
$X\ni x \mapsto \mathfrak{L}_{b_x}(g) \in \mathds{R}$, holding for
all $g\in C_{\rm b}^{+}(X)$.  The remaining assumptions are being
made to control the production of new particles, of which (ii) and
(iii) are related to the properties of $b_x(\Gamma^n)$, $n\in
\mathds{N}_0$, see (\ref{19F}). In general, (ii) and (iii) may be
quite independent as the choice of $\delta(x)$ leaves enough
possibilities to modify $n(x)$. However, in some cases, $\delta(x)$
and $n(x)$ can be expressed through each other. For instance, if
$b_x$ is a Poisson measure -- which is completely determined by its
first correlation measure that appears in (\ref{19F}), see, e.g.,
\cite[page 45]{Dawson} -- then $\delta(x) = e^{-n(x)}$. In this
case, (ii) follows by (iii) with $n_* = - \log \delta_*$. The role
of (iv) is to control the dispersal of offsprings, and thus the
nonlocality of the process. To illustrate its role, we take
$X=\mathds{R}$ and
\[
\bar{\beta}_x( d y) := \beta^{(1)}_x( d y) /n(x)= \frac{1}{2r}
\mathds{1}_{[x-r, x+r]} (y) d y, \qquad r>0.
\]
Then $\psi(y) = (1-\delta_*)e^{-\alpha |y|}$ satisfies
\[
\int_X \psi(y) \bar{\beta}_x( d y) \leq \left(\frac{e^{\alpha r} -
e^{-\alpha r}}{2 \alpha r} \right) \psi(x),
\]
which yields (\ref{Q}) with $m=\sinh(\alpha r)/\alpha r>1$. Note
 that this $m$ can be made arbitrarily close to one by
taking small enough either $r$ or $\alpha$. The former corresponds
to a short dispersal, whereas by choosing small $\alpha$ one makes
$\Gamma^\psi$ -- and hence $\mathcal{P}(\Gamma^\psi)$ -- smaller,
cf. (\ref{8Fa}) and (\ref{11F}).

\section{The Kolmogorov Equation}

Our aim in this section is to solve (\ref{KE}) and to prepare
solving (\ref{FPE}).

\subsection{Solving the log-Laplace equation}

 As is
typical for the theory of measure-valued branching processes, see
\cite{BL,BL1} and \cite[Chapt. 4]{Dawson}, the main point of their
constructing  is solving a nonlinear evolution equation, often
called `log-Laplace equation', see \cite[pp. 60, 61]{Dawson}. We
approach this by defining
\begin{equation}
  \label{20F}
  C_\psi(X) = \left\{\phi\in C_{\rm b}(X): \forall x\in X \quad 0<  c_\phi \psi(x) \leq 1-\phi(x)=:\theta(x) \leq 1- \delta(x)
  \right\},
\end{equation}
i.e., each $\theta=1-\phi$ has its own lower bound, whereas the
upper bound is one and the same for all such functions. Notably, by
item (iii) of Assumption \ref{Ass1} it follows that each $\phi\in
C_\psi(X)$ satisfies
\begin{equation}
  \label{20Fz}
 \phi(x)  \geq 1-\psi(x) \geq \delta_*.
\end{equation}
Let us prove that $(\Phi \phi)(x) \geq \delta(x)$, holding for each
$\phi\in C_\psi(X)$. Since $\phi(y) \geq 0$, by (\ref{21F}) we have
\begin{gather}
\label{Lvv}
 (\Phi \phi)(x) \geq \int_{\Gamma^0} b_x ( d \xi) = b_x(\Gamma^0)
 = \delta(x) \geq 1-\psi (x)\geq \delta_*,
\end{gather}
see item (iii) of Assumption \ref{Ass1}. Moreover, by (\ref{13F})
and (\ref{20Fz}) it follows that
\begin{equation}
  \label{20Fy}
g(x) \leq -\frac{1}{\psi(x)} \log (1-\psi(x)) = \sum_{n=1}^\infty
\frac{[\psi(x)]^{n-1}}{n} \leq -
\frac{\log(1-\delta_*)}{1-\delta_*}=: g_*.
\end{equation}
Both (\ref{20Fz}) and (\ref{20Fy}) holding for all $x\in X$.

 Now for
$T>0$, by $\mathcal{C}^T$ we denote the Banach space of continuous
maps $[0,T]\ni t\mapsto \varphi_t \in C_{\rm b}(X)$, equipped with
the norm
\begin{equation}
  \label{24F}
  \|\varphi\|_T = \sup_{t\in [0,T]} \sup_{x\in X}|\varphi_t(x)|.
\end{equation}
We also set
\[
\mathcal{C}_\psi^T =\{ \varphi\in \mathcal{C}^T: \varphi_t \in
C_\psi (X), \ t\in [0,T]\},
\]
and
\begin{equation}
  \label{d}
\mathcal{C}_\psi^T (\phi) =\{\varphi\in \mathcal{C}_\psi^T:
\varphi_0 = \phi, \ \  \varphi_t(x) \leq 1- c_\phi e^{-t} \psi(x)\},
\qquad \phi \in C_\psi (X),
\end{equation}
where $c_\phi$ is the same as in (\ref{20F}) for this $\phi$.
Clearly, $\mathcal{C}_\psi^T (\phi)$ is a closed subset of
$\mathcal{C}_\psi^T $. Indeed, let $\{\varphi_n\}_{n\in
\mathds{N}}\subset \mathcal{C}_\psi^T (\phi)$ be
$\|\cdot\|_T$-convergent to a certain $\varphi\in \mathcal{C}^T$.
Then $\varphi_0 =\phi$ and $\varphi_t(x) \leq 1- c_\phi e^{-t}
\psi(x)$ since $(\varphi_n)_t(x) \to \varphi_t(x)$ as $n\to
+\infty$, holding for all $t\in [0,T]$ and $x\in X$. Now we define
\begin{equation}
  \label{25F}
  (K \varphi)_t(x) = \varphi_0 (x) e^{-t} + \int_0^t e^{-(t-s)}(\Phi
  \varphi_s)( x) d s.
\end{equation}
\begin{proposition}
  \label{1Fpn}
 Let $n_*$ introduced in Assumption \ref{Ass1} and $T$ satisfy $n_* (1-e^{-T}) <1$. Then for each $\phi\in
 C_\psi(X)$, the map $K$ has a unique fixed point $\varphi\in
 \mathcal{C}_\psi^T(\phi)$.
 \end{proposition}
\begin{proof}
We begin by showing that $K:
\mathcal{C}_\psi^T(\phi)\to\mathcal{C}_\psi^T(\phi)$ for each $T>0$.
Clearly, $x\mapsto (K\varphi)_t(x)$ is continuous and  $(K\varphi)_0
=\phi$ whenever $\varphi \in \mathcal{C}_\psi^T(\phi)$. The
continuity of $t\mapsto \Phi \varphi_t$ follows by the estimate, see
(\ref{21F}),
\begin{eqnarray}
& & \label{26F} \left|(\Phi \varphi_s)(x) - (\Phi \varphi_u)(x)
\right|  \leq  \int_{\Gamma_0} \left| \prod_{y\in \xi} \varphi_s(y)
-
\prod_{y\in \xi} \varphi_u(y) \right| b_x( d\xi) \\[.2cm] \nonumber
& & \quad \leq  \sup_{y\in X} |\varphi_s(y) - \varphi_u (y)|
\int_{\Gamma_0} |\xi| b_x(d\xi) \leq n_* \sup_{y\in X} |\varphi_s(y)
- \varphi_u (y)|.
\end{eqnarray}
This also yields the continuity of $t \mapsto (K \varphi)_t$. In
obtaining (\ref{26F}) we have used the following evident estimate
\[
|a_1 a_2 \cdots a_n - b_1 b_2 \cdots b_n | \leq n \max_{i} |a_i -
b_i|, \qquad a_i, b_i \in [0,1].
\]
Furthermore,
\[
0< (K\varphi)_t(x) \leq \phi(x) e^{-t} + (1-e^{-t}) = 1 -
(1-\phi(x))e^{-t} \leq 1
\]
which yields
\begin{equation}
  \label{20FB}
  1 - (K\varphi)_t(x) \geq e^{-t} \theta(x) \geq e^{-t} c_\phi \psi(x) =:
  c_\phi(t)\psi(x),
\end{equation}
and hence the validity of the upper estimate assumed in (\ref{d}).
Similarly as in  (\ref{Lvv}) we have
\begin{eqnarray*}
(\Phi \varphi_s) (x) \geq  b_x(\Gamma^0) =\delta(x) \geq 1 -
\psi(x),
\end{eqnarray*}
where we used also item (iii) of Assumptions \ref{Ass1}. By means of
this estimate applied in (\ref{25F}) we then get
\begin{gather*}
  (K\varphi)_t(x) \geq \phi(x) e^{-t} + (1-e^{-t} )\delta (x)  \\[.2cm] \nonumber \geq
(1-\psi(x) ) + e^{-t} (\phi(x) - \delta(x)) \geq 1-\psi(x),
\end{gather*}
as $\phi \in C_\psi(X)$. Thus, $K:
\mathcal{C}_\psi^T(\phi)\to\mathcal{C}_\psi^T(\phi)$. Let us show
that it is a contraction. To this end, similarly as in (\ref{26F})
we obtain, see also (\ref{24F}),
\[
\|K\varphi - K \tilde{\varphi}\|_T \leq n_* (1-e^{-T}) \|\varphi -
\tilde{\varphi}\|_T, \qquad \varphi, \tilde{\varphi} \in
\mathcal{C}_\psi^T (\phi).
\]
Now the proof follows by Banach's contraction principle.
\end{proof}
Next we consider the following nonlinear equation
\begin{equation}
  \label{29F}
\frac{\partial}{\partial t} \phi_t (x) = - \phi_t (x) + (\Phi
\phi_t)(x), \qquad \phi_0 = \phi.
\end{equation}
It is a nonlocal analog of the log-Laplace equation -- a standard
object in the theory of branching processes, see, e.g., \cite[page
61]{Dawson}. By a classical solution of (\ref{29F}) we will
understand a map $\mathds{R}_{+} \ni t \mapsto \phi_t \in C_{\rm
b}(X)$ which is everywhere continuously differentiable and satisfies
both equalities mentioned therein.
\begin{lemma}
  \label{LUlm}
For each $\phi\in C_\psi(X)$, (\ref{29F}) has a unique solution
$t\mapsto \phi_t \in C_\psi(X)$ which satisfies
\begin{equation}
  \label{bounds}
  c_\phi (t)\psi(x)\leq  1 - \phi_t(x) =:\theta_t(x) \leq \psi(x),
\end{equation}
with $c_\phi(t)$ defined in (\ref{20FB}). For $n_* <1$, this
solution tends to $\phi_\infty (x) \equiv 1$ as $t\to +\infty$ in
the norm of $C_{\rm b}(X)$.
\end{lemma}
\begin{proof}
We begin by fixing $T>0$ such that the contraction condition
$n_*(1-e^{-T}) < 1$ is satisfied. Then integrating in (\ref{29F}) we
arrive at the following integral equation
\begin{equation}
  \label{30F}
\phi_t (x) = \phi(x)e^{-t} + \int_0^t e^{-(t-s)} (\Phi\phi_s)(x) ds,
\end{equation}
the set of solutions of which on $[0,T]$ coincides with the set of
fixed points of $K:\mathcal{C}^T(\phi)\to \mathcal{C}^T(\phi)$
established in Proposition  \ref{1Fpn}. The continuous
differentiability of $t\mapsto \phi_t\in C_{\rm b}(X)$ follows by
continuity $s\mapsto \Phi\psi_s$, which in turn follows by
(\ref{26F}). Thus, each solution of (\ref{30F}) solves also
(\ref{29F}), which yields the existence of the solution in question
on the time interval $[0,T]$. For $n_*\leq 1$, the contraction
condition is satisfied for any $T>0$; hence, the aforementioned
solution is global in time. For $n_*>1$, we proceed as follows. For
$t_1+t_2 \leq T$, we rewrite (\ref{30F}) as follows
\begin{gather}
\label{31F} \phi_{t_1+t_2} (x) = e^{-t_2} \bigg{(} \phi(x) e^{-t_1}
+
\int_0^{t_1} e^{-(t_1 - s)} (\Phi \phi_s)(x) d s \bigg{)} \\[.2cm]
\nonumber + \int_{t_1}^{t_1+t_2} e^{-(t_2+t_1-s)} (\Phi \phi_{s})
(x) ds
\\[.2cm] \nonumber =
\phi_{t_1}(x) e^{-t_2} + \int_{0}^{t_2} e^{-(t_2-s)} (\Phi
\phi_{t_1+ s}) (x) ds.
\end{gather}
Since the contraction condition is independent of the initial
condition in (\ref{29F}), by (\ref{31F}) one can continue the
solution obtained above to any $t>0$. Indeed, let $\phi_t$ be the
solution on $[0,T]$. Let also $\phi^1_t\in
\mathcal{C}^T_\psi(\phi^1)$ be the solution of (\ref{29F}) on the
same $[0,T]$ with the initial condition $\phi^1_t := \phi_{T/2}$. By
the uniqueness established in Lemma \ref{LUlm} it follows that these
 two solutions satisfy $\phi_{t+T/2} = \phi^1_t$ for $t\in
[0,T/2]$. Hence, the function $\phi_t \mathds{1}_{[0,T/2]}(t) +
\phi^{1}_{t-T/2} \mathds{1}_{[T/2, 3T/2]}(t)= \phi_t
\mathds{1}_{[0,T]}(t) + \phi^{1}_{t-T} \mathds{1}_{[T, 3T/2]}(t)$ is
the unique solution of (\ref{30F}) (hence of (\ref{29F})) on $[0,
3T/2]$. The further continuation goes in analogous way.

For $n_*<1$, we define $\vartheta_s = e^s \|1 - \phi_s\| = e^s
\sup_{x\in X} (1 - \phi_s(x))$. By (\ref{30F}) we then get
\begin{equation*}
\vartheta_t \leq \vartheta_0 + n_*  \int_0^t \vartheta_s ds.
\end{equation*}
which by Gr\"onwall's inequality yields,
\[
\|1-\phi_t\| \leq \|1-\phi\|e^{- (1-n_*)t},
\]
and thereby the convergence in question. Note that $\phi_\infty$
does not belong to $C_\psi(X)$ as it fails to obey the upper bound
$\phi(x) \leq 1 - c_\phi \psi(x)$ with $c_\phi>0$, see (\ref{20F}).
However, it belongs to the closure of this set, and is a stationary
solution of (\ref{29F}).
\end{proof}
\begin{remark}
  \label{qrk}
 By (\ref{31F}) it follows that the solution of (\ref{29F}) -- which is
 a  nonlinear Cauchy problem in the Banach space $C_{\rm b}(X)$ --
 is given by a continuous semigroup of nonlinear operators, say
 $\{\rho_t\}_{t\geq 0}$, in the form $\phi_t = \rho_t(\phi_0)$, $\phi_t
 \in C_\psi(X)$. If one writes $\phi_t
 \in C_\psi(X)$ in the form $\phi_t (x) = \exp( - g_t (x)
 \psi((x))$, see (\ref{13F}), then the map $g\mapsto g_t$ also has
 the flow property. It defines a continuous
 semigroup of nonlinear operators $\{r_t\}_{t\geq 0}$ such that $g_t
 = r_t(g_0)$. It is known as the log-Laplace semigroup, see
 \cite[page 60]{Dawson}.
\end{remark}
We conclude this subsection by establishing the following useful
properties of the solution $\phi_t$.
\begin{lemma}
  \label{LU1lm}
Let $\phi_t= 1 - \theta_t -e^{-g_t\psi}$ be the solution as in Lemma
\ref{LUlm}. Then, for each $t\geq 0$, $u>0$ and all $x\in X$, the
following holds
\begin{eqnarray}
  \label{Lv}
  & (a) & \quad  |\phi_{t+u}(x) - \phi_t (x)| = |\theta_{t+u}(x) - \theta_t (x)| \leq 2 u
 \psi(x), \\[.2cm] \nonumber
& (b) & \quad |g_{t+u}(x) - g_t (x)| \leq 2 u/\delta_*, \\[.2cm] \nonumber
& (c) & \quad |(\Phi\phi_{t+u})(x) - (\Phi\phi_t) (x)| \leq  2 u n_*
m \psi(x).
\end{eqnarray}
\end{lemma}
\begin{proof}
 By (\ref{29F}) we have
\begin{gather}
  \label{Lv1}
|\phi_{t+u}(x) - \phi_t(x)|\leq \int_0^u | \phi_{t+s}(x) - (\Phi
\phi_{t+s}) (x) | d s \\[.2cm] \nonumber= \int_0^u | \theta_{t+s}(x) - (1-
(\Phi\phi_{t+s})(x))| d s \leq 2 \psi (x) u,
\end{gather}
where we have used (\ref{bounds}) and (\ref{Lvv}). To prove (b), we
denote $$h^+(x) = \max\{g_{t+u}(x) \psi(x); g_{t}(x) \psi(x)\},
\quad h^-(x) = \min\{g_{t+u}(x) \psi(x); g_{t}(x) \psi(x)\}.$$ Then,
cf. (\ref{13F}),
\begin{eqnarray*}
& & |\phi_{t+u}(x) - \phi_t(x)|  =  e^{-h^{+}(x)} \left[
e^{h^{+}(x)-h^{-}(x)}-1\right] \\[.2cm] \nonumber  & & \qquad \quad \geq e^{-h^{+}(x)}|g_{t+u}(x) -
g_t (x)| \psi (x)  \\[.2cm] \nonumber  & & \qquad \quad \geq
\min\{\phi_{t+u}(x) ; \phi_t(x)\} |g_{t+u}(x) - g_t (x)| \psi (x) ,
\end{eqnarray*}
which yields case (b) of (\ref{Lv}) by (\ref{Lv1}) and (\ref{20Fz}).
Next, similarly as in (\ref{26F}) we get
\begin{gather*}
\left|(\Phi\phi_{t+u})(x) - (\Phi\phi_{t})(x) \right|
 \leq \int_{\Gamma_0} \left(\sum_{y\in \xi} \left|\phi_{t+u}(y) - \phi_{t}(y)
\right| \right) b_x(d\xi) \\[.2cm] \nonumber \leq  2 u \int_X
\psi(y) \beta_x^{(1)} (d y) \leq 2 u n_* m \psi(x),
\end{gather*}
where we used (\ref{Lv1}) and (\ref{Q}), see also item (i) of
Assumption \ref{Ass1}.
\end{proof}

\subsection{Basic estimates}
In defining $L$, we employ a number of estimates which we derive
now. For $\phi=e^{-g \psi}\in C_\psi(X)$ we set, see (\ref{13F}),
\begin{equation}
  \label{15F}
  F^\phi(\gamma) = \prod_{x\in \gamma} \phi(x) = \exp\left( -\sum_{x\in \gamma} g(x)
  \psi(x)\right) = G^g (\gamma),
\end{equation}
where $G^g (\gamma)$ is as in (\ref{N6}).
\begin{proposition}
  \label{2Fpn}
Let $F^\phi$ be as in (\ref{15F}) with $\phi\in C_\psi(X)$, see
(\ref{20F}). Then, for each $\gamma\in \Gamma^\psi$, the following
holds
\begin{equation}
  \label{esti}
\left|LF^\phi(\gamma)\right| \leq \frac{2}{e\delta_* c_\phi},
\end{equation}
where $c_\phi$ defines the lower bound in (\ref{20F}). By
(\ref{esti}) it then follows that $LF^\phi \in C_{\rm
b}(\Gamma^\psi)$.
\end{proposition}
\begin{proof}
By (\ref{6F}), and then by (\ref{20Fz}), (\ref{Lvv}) and
(\ref{20F}), we have
\begin{gather}
  \label{28F}
  |LF^\phi (\gamma)| \leq \sum_{x\in \gamma} F^\phi (\gamma\setminus
  x) \left| (\Phi \phi)(x) - \phi(x)\right| \\[.2cm] \nonumber \leq
  (F^\phi (\gamma)/\delta_*) \sum_{x\in \gamma}\bigg{(} \left| 1-(\Phi \phi)(x) \right| + \left| 1 -
  \phi(x)\right|\bigg{)}\\[.2cm] \nonumber \leq 2
  \Psi(\gamma) F^\phi (\gamma)/\delta_* \leq 2
  F^\phi(\gamma) e^{c_\phi \Psi(\gamma)} /(e\delta_*c_\phi) \leq 2/(e\delta_*c_\phi),
  \end{gather}
where $\Psi$ is as in (\ref{7F}). To get the latter two estimates in
(\ref{28F}), we proceeded as  follows. The first one was obtained
with the help of the estimate $\alpha \leq e^{\alpha -1}$, $\alpha
\geq 0$. Afterwards, we estimated
\begin{equation}
  \label{July13}
F^\phi (\gamma) \exp(c_\phi \Psi (\gamma)) = \prod_{x\in \gamma}
(1-\theta(x)) e^{c_\phi \psi(x)} \leq  \prod_{x\in \gamma} (1-c_\phi
\psi(x)) e^{c_\phi \psi(x)} \leq 1,
\end{equation}
see (\ref{20F}), which was used in the final step.  The continuity
of the map $\gamma \mapsto LF^\phi(\gamma)$ follows by the very
definition of the topology of $\Gamma^\psi$.
\end{proof}
As in (\ref{20F}) we do not restrict the lower bounds,  the
right-hand side of (\ref{esti}) can be arbitrarily large for
 small enough $c_\phi$.
\begin{corollary}
  \label{Julyco}
For a given $\phi \in C_\psi(X)$, let $\phi_t$ be the solution
mentioned in Lemma \ref{LUlm}. Then for each $\gamma\in
\Gamma^\psi$, the map $t\mapsto F^{\phi_t}(\gamma)$ is continuously
differentiable on $\mathds{R}_{+}$ and the following holds
\begin{equation}
  \label{July}
  \frac{d}{dt} F^{\phi_t}(\gamma) = LF^{\phi_t}(\gamma).
\end{equation}
\end{corollary}
\begin{proof}
Recall that $\phi_t(x) = \exp\left(- g_t(x) \psi(x) \right)$. Then
the continuous differentiability of $t \mapsto g_t(x)$ follows by
the analogous property of $t \mapsto \phi_t(x)$, see Lemma
\ref{LUlm}. Indeed, by case (b) of (\ref{Lv})  it follows that the
derivative of the former map is bounded uniformly in $x$, which
yields that the map $t \mapsto \sum_{x\in \gamma} g_t(x) \psi(x)$ is
continuously differentiable for each $\gamma\in \Gamma^\psi$. This
implies the same property for the map $t\mapsto F^{\phi_t}(\gamma)$,
as well as the boundedness of the map $\gamma \mapsto (d/dt)
F^{\phi_t}(\gamma)$. The latter is proved analogously as in
(\ref{28F}). Finally, the equality in (\ref{July}) follows by the
fact that $\phi_t$ solves (\ref{29F}).
\end{proof}
Our next step is obtaining a number of useful estimates
characterizing the map $t \mapsto F^{\phi_t}(\gamma)$.
 \begin{lemma}
   \label{Lv1lm}
 For a given $\phi\in C_\psi(X)$, let $\phi_t$ be the solution of (\ref{29F}), see Lemma \ref{LUlm}.
 Then, for each $t\geq 0$, $u>0$ and $\gamma\in \Gamma^\psi$, the
 following holds
\begin{equation*}
\left|F^{\phi_{t+u}}(\gamma) - F^{\phi_{t}}(\gamma)  \right| \leq
\frac{2ue^{t+u}}{e\delta_* c_\phi}.
\end{equation*}
 \end{lemma}
\begin{proof}
We fix $t$ and $u$ and define
\begin{gather*}
  H_s (\gamma) = \sum_{x\in \gamma} g_s(x) \psi(x), \quad  H^{+}
  (\gamma) = \max\{ H_{t+u} (\gamma);  H_t (\gamma)\}, \\[.2cm] \nonumber H^{-}
  (\gamma) = \min\{ H_{t+u} (\gamma);  H_t (\gamma)\}.
\end{gather*}
Then
\begin{eqnarray}
  \label{Lv14}
& & \left|F^{\phi_{t+u}}(\gamma) - F^{\phi_{t}}(\gamma)  \right| =
e^{-H^{+}(\gamma)} \left[e^{H^{+}(\gamma)-H^{-}(\gamma)} -1
\right]\\[.2cm] \nonumber & &  \leq \max\{F^{\phi_{t+u}}(\gamma);
F^{\phi_{t}}(\gamma)\} \sum_{x\in \gamma} |g_{t+u}(x) -
g_t(x)|\psi(x) \\[.2cm] \nonumber & &  \leq \frac{2u}{\delta_*}\Psi(\gamma)\prod_{x\in \gamma} ( 1 -
c_\phi(t+u) \psi(x)) \\[.2cm] \nonumber &  &\leq \frac{2u}{e\delta_* c_\phi(t+u)}
\prod_{x\in \gamma} ( 1 - c_\phi(t+u) \psi(x))e^{c_\phi(t+u)\psi(x)}\\[.2cm] \nonumber &
&\leq \frac{2ue^{t+u}}{e\delta_* c_\phi},
\end{eqnarray}
which completes the proof, see (\ref{Lv}), (\ref{bounds}) and
(\ref{28F}).
\end{proof}
\begin{lemma}
  \label{Lv2lm}
Let $\phi$, $t$ and $u$ be as in Lemma \ref{Lv1lm}. Then there
exists $C_\phi>0$ such that, for all $\gamma\in \Gamma^\psi$, the
following holds
\begin{equation}
  \label{Lv15}
\left|(L F^{\phi_{t+u}})(\gamma) - (L F^{\phi_{t}})(\gamma)  \right|
\leq C_\phi u e^{2(t+u)}.
\end{equation}
\end{lemma}
\begin{proof}
As in (\ref{28F}), for fixed $t$ and $u$ we have
\begin{eqnarray}
  \label{Lv16}
& & \left| (L F^{\phi_{t+u}})(\gamma) - (L
F^{\phi_{t}})(\gamma)\right| \leq K_1(\gamma) + K_2(\gamma) +
K_3(\gamma),\\[.2cm] \nonumber & & K_1(\gamma) := \sum_{x\in \gamma} \left|F^{\phi_{t+u}}(\gamma\setminus x) -
F^{\phi_{t}}(\gamma\setminus x) \right| \left|(\Phi \phi_{t+u})(x)-
\phi_{t+u}(x)\right|, \\[.2cm] \nonumber & & K_2(\gamma) := \sum_{x\in \gamma}F^{\phi_{t}}(\gamma\setminus x)
\left| (\Phi \phi_{t+u})(x)- (\Phi \phi_{t})(x)\right|, \\[.2cm] \nonumber & & K_3(\gamma) := \sum_{x\in \gamma}F^{\phi_{t}}(\gamma\setminus x)
\left|  \phi_{t+u}(x)-  \phi_{t}(x)\right|.
\end{eqnarray}
By (\ref{Lvv}) and (\ref{20FB}) we have
\[
\frac{1}{1-c_\phi(t+u) \psi(x)} \leq \frac{1}{1-c_\phi \psi(x)} \leq
\frac{1}{1- \psi(x)} \leq \frac{1}{\delta_*}.
\]
Then proceeding as in obtaining the second inequality in
(\ref{Lv14}), we arrive at
\begin{eqnarray}
  \label{Lv17}
& & \left|F^{\phi_{t+u}}(\gamma\setminus x) -
F^{\phi_{t}}(\gamma\setminus x) \right| \leq \frac{2u}{\delta_*}
\Psi(\gamma\setminus x) \prod_{y\in \gamma\setminus
x}\left(1-c_\phi(t+u) \psi(y) \right)\\[.2cm] \nonumber & & \qquad  \leq
\frac{2u}{\delta^2_*} \Psi(\gamma) \prod_{y\in
\gamma}\left(1-c_\phi(t+u) \psi(y) \right)
\end{eqnarray}
Next, by (\ref{20Fz}) and (\ref{Lvv}) we have
\[
\left|(\Phi \phi_{t+u})(x) - \phi_{t+u}(x) \right| \leq
\left|1-(\Phi \phi_{t+u})(x) \right| + \left|1- \phi_{t+u}(x)
\right| \leq 2 \psi(x).
\]
We use the latter estimate and (\ref{Lv17}) to obtain
\begin{eqnarray}
  \label{Lv18}
  K_1(\gamma) & \leq & \frac{4u}{\delta^2_*} \Psi^2(\gamma) \prod_{y\in
\gamma}\left(1-c_\phi(t+u) \psi(y) \right) \\[.2cm] \nonumber &
\leq& \frac{16u}{(e \delta_* c_\phi(t+u))^2}\prod_{y\in
\gamma}\left(1-c_\phi(t+u) \psi(y) \right) e^{c_\phi(t+u) \psi(y)} \\[.2cm] \nonumber &
\leq& \frac{16u}{(e \delta_* c_\phi)^2} e^{2(t+u)}.
\end{eqnarray}
By (\ref{Lv}) we have
\begin{eqnarray}
  \label{Lv19}
K_2(\gamma) & \leq & \frac{1}{\delta_*} F^{\phi_t}(\gamma)
\sum_{x\in \gamma} \left|(\Phi\phi_{t+u})(x) - (\Phi\phi_{t})(x)
 \right| \\[.2cm] \nonumber & \leq & \frac{2un_*m}{\delta_*}
 \Psi(\gamma)F^{\phi_t}(\gamma) \leq \frac{2un_*m}{e\delta_*
 c_\phi}e^t.
\end{eqnarray}
Similarly,
\begin{equation}
  \label{Lv20}
K_3(\gamma)  \leq \frac{1}{\delta_*} F^{\phi_t}(\gamma) \sum_{x\in
\gamma} \left|\phi_{t+u}(x) - \phi_{t}(x)
 \right| \leq \frac{2u}{e\delta_*
 c_\phi}e^t.
\end{equation}
Now we use (\ref{Lv18}), (\ref{Lv19}), (\ref{Lv20}) in (\ref{Lv16}),
and thus obtain (\ref{Lv15}) with $$C_\phi =
\frac{2u(n_*m+1)}{e\delta_* c_\phi}  + \frac{16u}{(e \delta_*
c_\phi)^2},$$ which completes the proof.
\end{proof}

\subsection{Solving the Kolmogorov equation}
Now we can turn to solving (\ref{KE}) and preparing to solving the
main equation (\ref{FPE}). Set
\begin{equation}
  \label{Lv21}
 E^0
(\Gamma^\psi):= {\rm l.s.} \{F^\phi: \phi \in C_\psi(X)\},
\end{equation}
where l.s. $=$ linear span. It is a subset of the Banach space
$C_{\rm b}(\Gamma^\psi)$ equipped with the norm
\[
\|F\|:=\sup_{\gamma\in \Gamma^\psi}|F(\gamma)| .
\]
\begin{remark}
  \label{ark}
The set $E^0 (\Gamma^\psi)$ has all the properties stated in
Proposition \ref{N2pn}. This follows by the fact that the family
$\{G^v : v \in \mathcal{V}\}$ mentioned therein is a subset of $E^0
(\Gamma^\psi)$, see (\ref{15F}). In particular, it contains constant
functions and each measurable and bounded $F:\Gamma^\psi \to
\mathds{R}$ can be obtained as the bp-limit of a sequence of the
elements of $E^0(\Gamma^\psi)$.
\end{remark}
By (\ref{esti}) we know that $L: E^0 (\Gamma^\psi) \to C_\psi(X)$.
In view f this, we introduce
\begin{equation}
  \label{Lv25}
\|F\|_L = \|F \| + \|LF\|, \qquad F\in E^0 (\Gamma^\psi),
\end{equation}
that is, $\|\cdot \|_L$ is the corresponding graph-norm. Thereby, we
define
\begin{equation}
  \label{Lv26}
  \mathcal{D}(L) = \overline{E^0 (\Gamma^\psi)}^L,
\end{equation}
i.e. $\mathcal{D}(L)$ is the closure of $E^0 (\Gamma^\psi)$ in the
graph-norm, and thus the operator $(L,\mathcal{D}(L))$ is closed.
Below -- in particular, in (\ref{KE}) -- by $L$ we will mean this
operator.

Following \cite[page 108]{Arendt} by a classical solution of
(\ref{KE}) we will understand a continuously differentiable map
$\mathds{R}_{+} \ni t \mapsto F_t  \in \mathcal{D}(L)\subset C_+{\rm
b}(\Gamma^\psi)$ such that both equalities in (\ref{KE}) are
satisfied.
\begin{lemma}
  \label{1tm}
  For each $\phi\in C_{\psi}(X)$, the map $t\mapsto F^{\phi_t}$ is a
  classical solution of the Cauchy problem in (\ref{KE}) with
  $F_0=F^\phi$. For $n_* <1$, this solution satisfies $F_t (\gamma)
  \to F_{\infty}(\gamma)$ as $t\to +\infty$, where $F_{\infty}(\gamma) \equiv 1$ and
  the convergence is to hold for all $\gamma\in \Gamma^\psi$.
\end{lemma}
\begin{proof}
We begin by noting that the map $t\mapsto F^{\phi_t}$ has the flow
property related to (\ref{31F}) and remarking that (\ref{esti})
implies
\begin{equation}
  \label{estia}
  \|L F^{\phi_t}\| \leq \frac{2}{e\delta_* c_\phi}e^t,
\end{equation}
see (\ref{bounds}) and (\ref{20FB}). By (\ref{Lv15}) we know that
the map $t \mapsto L F^{\phi_t} \in C_{\rm b} (\Gamma^\psi)$ is
continuous and hence Bochner-integrable on each interval
$[a,b]\subset \mathds{R}_{+}$. Therefore,
\[
\int_0^t L F^{\phi_u} d u \in C_{\rm b} (\Gamma^\psi),
\]
and its value at $\gamma\in \Gamma^\psi$ satisfies
\begin{equation}
  \label{estib}
\int_0^t L F^{\phi_u} (\gamma)  d u = F^{\phi_t}(\gamma) -
F^{\phi}(\gamma),
\end{equation}
that readily follows by Corollary \ref{Julyco}. In view of the
aforementioned flow property, the proof of the first part of this
statement can be done by showing that
\begin{equation}
  \label{July17}
  \kappa_\phi(t) := \left\|\frac{1}{t}\left( F^{\phi_t} -
F^{\phi} \right) - LF^{\phi} \right\| \to 0, \qquad t\to 0.
\end{equation}
By (\ref{estib}) and (\ref{Lv15}) for $t\in (0,1)$ we have
\begin{gather*}
 \kappa_\phi(t) = \sup_{\gamma \in \Gamma^\psi} \left| \frac{1}{t}\left( F^{\phi_t} (\gamma) -
F^{\phi} (\gamma) \right) - LF^{\phi}(\gamma)\right| \\[.2cm]
\nonumber \leq  \frac{1}{t} \sup_{\gamma \in \Gamma^\psi}\int_0^t
\left| LF^{\phi_u}(\gamma)- LF^{\phi}(\gamma)\right|du \leq (e^2
C_\phi/2) t,
\end{gather*}
which yields (\ref{July17}). In the subcritical case $n_*<1$, the
stated convergence follows by the concluding part of Lemma
\ref{LUlm}.
\end{proof}
Since the map $t\mapsto F^{\phi_t} \in C_{\rm b}(\Gamma^\psi)$ is
continuous and bounded (by one), the Bochner integral
\begin{equation}
  \label{Lv22}
F^\phi_\lambda = \int_0^{+\infty} e^{-\lambda t} F^{\phi_t} d t,
\qquad \phi \in C_\psi(X),
\end{equation}
exists for all $\lambda >0$. In view of (\ref{Lv15}) and
(\ref{estia}), the map $t \mapsto LF^{\phi_t}\in C_{\rm
b}(\Gamma^\psi)$ is continuous and absolutely $e^{-\lambda t} d
t$-integrable for all $\lambda >1$. This observation leads us to the
following fact.
\begin{lemma}
  \label{2Flm}
For each $\phi\in C_\psi (X)$ and $\lambda>1$, it follows that
$F^\phi_\lambda$ defined in (\ref{Lv22}) lies in $\mathcal{D}(L)$
and the following holds
\begin{eqnarray}
  \label{33F}
L F^\phi_\lambda  = \int_0^{+\infty} e^{-\lambda t} LF^{\phi_t} dt=
- F^{\phi}
   + \lambda F^\phi_\lambda .
\end{eqnarray}
\end{lemma}
\begin{proof}
In view of the existence of the Laplace transforms just discussed,
the facts that $L$ is closed and $F^{\phi_t}$ solves the Cauchy
problem in (\ref{KE}), see Lemma \ref{1tm}, both stated properties
follow by a direct application of \cite[Theorem 3.1.3, page
109]{Arendt}.
\end{proof}
Now we define
\begin{equation}
  \label{Lv26a}
  \mathcal{F}^0(\Gamma^\psi) = {\rm l.s.} \{F^\phi_\lambda: \phi\in C_\psi (X), \ \lambda
  >1\}, \qquad \mathcal{F}(\Gamma^\psi)
  =\overline{\mathcal{F}^0(\Gamma^\psi)}^L,
\end{equation}
where $F^\phi_\lambda$ are defined  in (\ref{Lv22}) and the closure
is taken in the graph-norm (\ref{Lv25}).
\begin{lemma}
  \label{Lv3lm}
It follows that $\mathcal{F}(\Gamma^\psi)= \mathcal{D}(L)$. Thereby,
$\mathcal{F}^0(\Gamma^\psi)$ is a core of $\mathcal{D}(L)$.
\end{lemma}
\begin{proof}
We begin by showing that
\begin{equation}
  \label{estic}
  E^0(\Gamma^\psi) \subset \mathcal{F}(\Gamma^\psi),
\end{equation}
i.e., each $F^\phi$, $\phi \in C_\psi(X)$, can be obtained as the
$\|\cdot\|_L$-limit of a sequence of the elements of
$\mathcal{F}^0(\Gamma^\psi))$. In fact, we are going to show that
\begin{equation}
  \label{Lvu}
\|\lambda F^\phi_\lambda - F^\phi\|_L \to 0 , \quad {\rm as} \
\lambda \to +\infty.
\end{equation}
To this end, with the help of the first equality in (\ref{33F}) for
$\lambda >1$ we write
\begin{gather}
  \label{Lv27}
 \left|\lambda (L F^\phi_\lambda) (\gamma) - (L F^\phi)(\gamma)\right| =
 \left| \int_0^{+\infty} \left[ (L F^{\phi_t})(\gamma) - (LF^\phi) (\gamma)\right] e^{-\lambda t} \lambda d
 t\right|\\[.2cm] \nonumber \leq \int_0^{+\infty} \left|(L F^{\phi_{\epsilon s}})(\gamma) - (LF^\phi) (\gamma)
 \right|e^{-s} ds, \quad \epsilon = 1/\lambda.
\end{gather}
Now we use here (\ref{Lv15}) with $t=0$, $u= \epsilon s$ and obtain
for $\epsilon < 1/2$ the following estimate
\begin{gather}
  \label{Lv28}
{\rm LHS(\ref{Lv27})} \leq \epsilon C_\phi \int_0^{+\infty} s
e^{-s(1-2\epsilon)} d s = \frac{\epsilon}{(1-2\epsilon)^2} C_\phi
\to 0, \ \ {\rm as} \ \epsilon \to 0.
\end{gather}
Next, by (\ref{33F}) -- and then by (\ref{esti}) --  we get
\begin{gather}
  \label{Lv29}
  \| \lambda F^\phi_\lambda - F^\phi\|= \|LF^{\phi}_\lambda\| \leq
  \int_0^{+\infty} \|L F^{\phi_t}\|e^{-\lambda t} dt \\[.2cm] \nonumber\leq
  \frac{2}{e\delta_* c_\phi} \int_0^{+\infty} e^{-(\lambda -1) t} d
  t = \frac{1}{\lambda -1}\left( \frac{2}{e\delta_* c_\phi}\right),
\end{gather}
where we have used the fact that $c_{\phi_t} = c_\phi(t) = c_\phi
e^{-t}$, see (\ref{20FB}). Then (\ref{Lvu}) readily follows by
(\ref{Lv28}) and (\ref{Lv29}). Now by (\ref{Lv26}), (\ref{Lv26a})
and (\ref{estic}) we get $\mathcal{D}(L)\subset
\mathcal{F}(\Gamma^\psi)$. At the same time, by Lemma \ref{2Flm} it
follows that $\mathcal{F}^0(\Gamma^\psi)\subset \mathcal{D}(L)$,
which yields the opposite inclusion $\mathcal{D}(L)\supset
\mathcal{F}(\Gamma^\psi)$.
\end{proof}
For $\phi\in C_\psi(X)$ and $t\geq 0$,  $\phi_t$ is also in
$C_\psi(X)$, see Remark \ref{qrk}, and thus we can consider
\begin{equation*}
  F^{\phi_t}_\lambda = F^{\rho_t(\phi)}_\lambda= \int_0^{+\infty} F^{\phi_t+s} e^{-\lambda s} d
  s,
\end{equation*}
and apply here Lemma \ref{2Flm} to $F^{\phi_t}_\lambda $. Then we
obtain the following
\begin{corollary}
  \label{July20co}
For each $\phi\in C_\psi (X)$ and $\lambda >1$, the map $t \mapsto
F_t= F^{\phi_t}_\lambda \in \mathcal{D}(L) \subset C_{\rm
b}(\Gamma^\psi)$ is a classical solution of the Cauchy problem
(\ref{KE}) with $F_0 = F^\phi_\lambda$.
\end{corollary}
\begin{proof}
Set
\begin{gather}
  \label{July20}
  V_\lambda (t) = F^{\phi_t}_\lambda - F^{\phi}_\lambda - t L
  F^\phi _\lambda.
\end{gather}
Then we apply here repeatedly (\ref{33F}) and
\begin{gather}
  \label{July21}
 V_\lambda (t) =  \int_0^t \left( L F^{\phi_u}_\lambda -  L F^{\phi}_\lambda\right) d u
= - \int_0^t\left( F^{\phi_u} - F^{\phi}\right) d u + \lambda
\int_0^t \left( F^{\phi_u}_\lambda - F^{\phi}_\lambda\right) d u
\\[.2cm] \nonumber = - \int_0^t\left( F^{\phi_u} - F^{\phi}\right) d u
- \frac{\lambda t^2}{2} F^\phi + \lambda^2
 \int_0^t \int_0^u
F^{\phi_v}_\lambda  d u dv =:  V^{(1)}_\lambda (t)+ V^{(2)}_\lambda
(t) + V^{(3)}_\lambda (t).
\end{gather}
By Lemma \ref{Lv1lm} and we have that
\[
\|V^{(1)}_\lambda (t) \| \leq \frac{t^2}{\delta_*c_\phi}, \qquad
\|V^{(2)}_\lambda (t) \| \leq \frac{\lambda t^2}{2} .
\]
Here we recall that $F^{\phi}(\gamma) \leq 1$ for each $\phi\in
C_\psi(X)$ and $\gamma\in \Gamma^\psi$, see (\ref{July13}). By this
and (\ref{Lv22}) we also have
\[
\|V^{(3)}_\lambda (t) \| \leq \frac{\lambda t^2}{2} .
\]
We apply this estimates in (\ref{July21}) and (\ref{July20}) and
then obtain
\[
\|V_\lambda (t) \|/t \to 0, \qquad t\to 0,
\]
which completes the proof.
\end{proof}

\section{The Fokker-Planck Equation}

\subsection{Solving the Fokker-Planck equation}
 Now we may turn to the probabilistic part of the topic. Recall
that we use probability measures on $\Gamma^\psi$ as states of the
studied system of branching particles.
\begin{definition}
  \label{N4adf}
By a solution of the Fokker-Planck equation  (\ref{FPE}) we
understand a map $\mathds{R}_{+}\ni t \mapsto \mu_t \in
\mathcal{P}(\Gamma^\psi)$ possessing the following properties: (a)
for each $F\in B_{\rm b}(\Gamma^\psi)$, the map $\mathds{R}_{+}\ni t
\mapsto \mu_t (F) \in \mathds{R}$ is measurable; (b) the equality in
(\ref{FPE}) holds for all  $F\in \mathcal{D}(L)$, where the latter
is defined in (\ref{Lv26}).
\end{definition}
\begin{theorem}
  \label{2tm}
  For each $\mu_0 \in \mathcal{P}(\Gamma^\psi)$, the Fokker-Planck
  equation (\ref{FPE}) has a unique solution in the sense of the
  definition given above. Moreover, this solution is weakly
  continuous, i.e., $\mu_t \Rightarrow \mu_s$ as $t\to s\in \mathds{R}_{+}$. In the subcritical case $n_*<1$,
$\mu_t \Rightarrow \mu_{\infty}$ as $t\to
  +\infty$, where $\mu_\infty$ is the measure supported on the singleton subset of $\Gamma^\psi$ consisting of the empty configuration, i.e.,
$\mu_\infty
  (\Gamma^0) =1$
\end{theorem}
The proof of this theorem is based, in particular, on the following
fact.
\begin{lemma}
  \label{alm}
Let a map $t \mapsto \mu_t$ satisfy condition (b) of Definition
\ref{N4adf}. Then it also satisfies (a), and hence is a solution of
(\ref{FPE}).
\end{lemma}
The proof of this statement in turn is based on the following
result, which has its own value.
\begin{proposition}
  \label{Au7Fpn}
Let $t \mapsto \mu_t \in \mathcal{P}(\Gamma^\psi)$ satisfy
(\ref{FPE}) for all $t_1,t_2$ and $F\in \mathcal{D}(L)$. Then, for
each $F\in \mathcal{F}^0(\Gamma^\psi)$, the map $t \mapsto \mu_t
(F)\in \mathds{R}$ is Lipschitz-continuous. The same is true also
for $F\in E^0 (\Gamma^\psi)$, see (\ref{Lv21}).
\end{proposition}
\begin{proof}
First, we rewrite (\ref{FPE}) in the form
\begin{equation}
  \label{FPEa}
\mu_{t_2} (F) = \mu_{t_1} (F) + \int_{t_1}^{t_2} \mu_s ( LF) d s,
\qquad 0\leq t_1 < t_2.
\end{equation}
Take $F=F^\phi_\lambda$, $\phi \in C_\psi(X)$, $\lambda >1$.  By
(\ref{15F}), and then by (\ref{Lv22}) and (\ref{33F}), we have $\|L
F^\phi_\lambda\|\leq 2$. Then by (\ref{FPEa}) one obtains
\[
|\mu_{t_2} (F^\phi_\lambda) - \mu_{t_1} ( F^\phi_\lambda)| \leq 2
|t_2 - t_1|.
\]
For $F = \sum_{n} \alpha_n F^{\phi_n}_{\lambda_n} \in
\mathcal{D}^0(L)$, this yields
\[
|\mu_{t_2} (F) - \mu_{t_1} ( F) | \leq 2\left(\sum_{n}|\alpha_n|
\right)|t_2 - t_1|.
\]
Now for $F=F^\phi$, $\phi\in C_\psi(X)$, see (\ref{Lv21}),
(\ref{Lv26}), by (\ref{esti}) we have
\[
|\mu_{t_2} (F^\phi) - \mu_{t_1} ( F^\phi) | \leq \frac{2}{e\delta_*
c_\phi}|t_2 - t_1|.
\]
The extension of the latter to the linear combinations of
$F^{\phi_n}$ can be done similarly as above.
\end{proof}
\noindent {\it Proof of Lemma \ref{alm}.} By Remark \ref{ark} we
know that $E^0(\Gamma^\psi)$ is $bp$-dense in $B_{\rm
b}(\Gamma^\psi)$. Then the measurability of $t\mapsto \mu_t(F)$,
$F\in B_{\rm b}(\Gamma^\psi)$ follows by the continuity (hence,
measurability) just proved. \hfill $\square$

\noindent {\it Proof of Theorem \ref{2tm}.} In view of the lemma
just proved, it remains to establish the existence and uniqueness of
solutions of (\ref{FPE}) with $F\in \mathcal{D}(L)$. First we prove
existence. Let $F$ be in $\mathcal{F}^0(\Gamma^\psi)$ which is the
core of $\mathcal{D}(L)$, see (\ref{Lv26a}). Since (\ref{FPE}) is
linear, it is enough to take $F= F^\phi_\lambda$ with $\phi\in
C_\psi(X)$ and $\lambda
>1$. For a given $\mu \in \mathcal{P}(\Gamma^\psi)$ and $t>0$, we
then set
\begin{equation}
  \label{37F}
  \mu_t (F^\phi_\lambda) = \mu(F_\lambda^{\phi_t}),
\end{equation}
which has to hold for all $\phi\in C_\psi(X)$ and $\lambda
>1$. Let us multiply both sides of (\ref{37F}) by $\lambda$ and then
pass to the limit $\lambda \to +\infty$. By (\ref{Lvu}) we then get
that $\mu_t (F^{\phi}) = \mu(F^{\phi_t})$, which uniquely determines
$\mu_t$ in view of property (iii) of Proposition \ref{N2pn} and the
uniqueness stated in Lemma \ref{LUlm}, see also Remark \ref{ark}.
Note that $\mu_t (F^\phi) \in (0,1)$ for all $t>0$.  Then by
(\ref{33F}) we have
\begin{gather}
  \label{38F}
  \int_{t_1}^{t_2} \mu_s ( LF^\phi_\lambda) d s = -  \int_{t_1}^{t_2}
  \mu_s (F^\phi) d s + \int_{t_1}^{t_2} \mu_s (\lambda  F^\phi_\lambda) d s
\\[.2cm] \nonumber =
-  \int_{t_1}^{t_2}
  \mu_s (F^\phi) d s + \int_{t_1}^{t_2} \int_0^{+\infty} \lambda e^{-\lambda
  t} \mu_s (F^{\phi_t}) d s d t,
  \end{gather}
where we used also Fubini's theorem.  Thereafter, by (\ref{37F}) and
the flow property, see (\ref{31F}), we get $\mu_s (F^{\phi_t}) =
\mu_{s+t} (F^\phi)$ and then use this in the second summand (name it
$\Upsilon$) of the last line of (\ref{38F}), then integrate by parts
and by Fubini's theorem obtain
\begin{gather*}
\Upsilon = \int_{t_1}^{t_2} \mu_s (F^\phi) d s + \int_{t_1}^{t_2}
\frac{d}{ds}\left( \int_0^{+\infty}e^{-\lambda t} \mu_{s+t} (F^\phi)
d t \right) d s \\[.2cm] = \int_{t_1}^{t_2} \mu_s (F^\phi) d s + \int_{t_1}^{t_2}
\frac{d}{ds} \mu_s \left( \int_0^{+\infty}e^{-\lambda t} F^{\phi_t}
d t \right) d s \\[.2cm] = \int_{t_1}^{t_2} \mu_s (F^\phi) d s +  \int_{t_1}^{t_2}
\frac{d}{ds} \mu_s (F^\phi_\lambda) d s \\[.2cm] = \int_{t_1}^{t_2} \mu_s (F^\phi) d s +
\mu_{t_2}(F^\phi_\lambda) - \mu_{t_1}(F^\phi_\lambda).
\end{gather*}
Now we plug this in (\ref{38F}) and get that the map $t\mapsto
\mu_t(F)$, $F\in \mathcal{F}^0(\Gamma^\psi)$, solves (\ref{FPEa}).
For $F\in \mathcal{D}(L)$, let $\{F_n\}_{n\in \mathds{N}}\subset
\mathcal{F}^0(\Gamma^\psi)$ be such that $\|F - F_n\|_L \to 0$ as
$n\to +\infty$. Then
\begin{gather*}
\left|\mu_{t_2} (F) - \mu_{t_1}(F) - \int_{t_1}^{t_2} \mu_s (LF) d s
\right| \leq \left|\mu_{t_2} (F-F_n) \right| + \left|\mu_{t_1}
(F-F_n) \right|\\[.2cm] \nonumber  + \int_{t_1}^{t_2} \left|\mu_s (LF - LF_n) \right|
d s \leq (t_2 - t_1 +2) \|F - F_n\|_L,
\end{gather*}
which yields that $t\mapsto \mu_t(F)$, $F\in \mathcal{D}(L)$ also
solves (\ref{FPEa}).

Assume now that $t \mapsto \tilde{\mu}_t$ is another solution of
(\ref{FPE}), and hence of (\ref{FPEa}), satisfying
$\tilde{\mu}_t|_{t=0} = \mu$. By Proposition \ref{Au7Fpn} the map $t
\mapsto \tilde{\mu}(F)$, $F\in \mathcal{F}^0(\Gamma^\psi)$ is
Lipschitz-continuous.  Then, for each $\lambda>1$ and $\phi\in
C_\psi(X)$, we have
\[
d \tilde{\mu}_s (F^\phi_\lambda) =\tilde{ \mu}_s (L F^\phi_\lambda)
d s,
\]
holding for Lebesgue-almost all $s\geq 0$.  Then
\begin{gather*}
  - \lambda \int_0^t e^{-\lambda s} \tilde{\mu}_s (F^\phi_\lambda) d s =
  \int_0^t  \tilde{\mu}_s (F^\phi_\lambda) d e^{-\lambda s} \\[.2cm] = \tilde{\mu}_t
  (F^\phi_\lambda) e^{-\lambda t} - \tilde{\mu}_0 (F^\phi_\lambda) - \int_0^t
  e^{-\lambda s} \tilde{\mu}_s (L F^\phi_\lambda) ds \\[.2cm] = \tilde{\mu}_t
  (F^\phi_\lambda)e^{-\lambda t} - \tilde{\mu}_0 (F^\phi_\lambda) - \lambda \int_0^t
  e^{-\lambda s} \tilde{\mu}_s (F^\phi_\lambda) d s + \int_0^t
  e^{-\lambda s} \tilde{\mu}_s (F^\phi) d s.
\end{gather*}
This yields
\[
\mu (F^\phi_\lambda) = \tilde{\mu}_0 (F^\phi_\lambda)= \tilde{\mu}_t
  (F^\phi_\lambda)e^{-\lambda t} + \int_0^t
  e^{-\lambda s} \tilde{\mu}_s (F^\phi) d s, \qquad \lambda >1,
\]
which after passing to the limit $t\to +\infty$ leads to
\begin{equation}
  \label{39F}
\mu (F^\phi_\lambda) = \int_0^{+\infty}
  e^{-\lambda s} \tilde{\mu}_s (F^\phi) d s,
\end{equation}
that holds for all $\lambda >1$. By the very definition in
(\ref{37F}) the map $t \mapsto \mu_t(F^\phi)$ is continuous; the
continuity of $t \mapsto \tilde{\mu}_t(F^\phi)$ was established in
Proposition \ref{Au7Fpn}. Both maps are bounded. By (\ref{Lv22}) and
(\ref{37F}), and then by (\ref{39F}), the Laplace transforms of both
these maps coincide. Therefore, by Lerch's theorem
$\mu_t(F^\phi)=\tilde{\mu}_t(F^\phi)$ for all $t>0$ and $\phi\in
C_\psi(X)$. As mentioned above, see Proposition \ref{N2pn}, the
class of functions $\{F^\phi:\phi \in C_\psi(X)\}$ is separating,
that means $\mu_t=\tilde{\mu}_t$, $t>0$ and hence the stated
uniqueness. The proof the weak convergence $\mu_t \Rightarrow \mu_s$
follows by (\ref{37F}) and the fact that $\{F^\phi:\phi \in
C_\psi(X)\}$ is also convergence determining, see again Proposition
\ref{N2pn}. It remains to prove that $\mu_t \Rightarrow \mu_\infty$
as $t\to +\infty$. Since the set $\{F^\phi:\phi \in C_\psi(X)\}$ is
convergence determining, to this end it is enough to show that
$\mu_t(F^\phi)\to \mu_\infty (F^\phi)=1$, holding for all $\phi \in
C_\psi(X)$. By (\ref{37F}) and the concluding statement of Lemma
\ref{1tm} we have
\[
\lim_{t\to +\infty} \mu_t(F^\phi) = \lim_{t\to +\infty}
\mu(F^{\phi_t}) = \mu(F_\infty)=1,
\]
which completes the whole proof. \hfill $\square$

\subsection{Concluding comments}
As mentioned above, our main aim in this work is to find a way of
describing branching in infinite particle systems. That is why we
restrict ourselves to the results stated in Theorem \ref{2tm}. A
direct consequence of Theorem \ref{2tm} is the existence of a Markov
process with values in $\Gamma^\psi$, that may be constructed by
means of the Markov transition function $p^t_\gamma$, see
\cite[pages 156, 157]{EK}, determined by its values on
$\{F^\phi:\phi\in C_\psi(X)\}$, cf. Remark \ref{ark}. These values
are given by the following formula
\[
p^\gamma_t(F^\theta) = F^{\phi_t}(\gamma), \qquad \gamma \in
\Gamma^\psi.
\]
It definitely has the branching property, cf. \cite[page 29]{Li},
\[
F^{\phi_t}(\gamma_1\cup \gamma_2) = F^{\phi_t}(\gamma_1)
F^{\phi_t}(\gamma_2), \qquad \gamma_1, \gamma_2 \in \Gamma^\psi.
\]
Then in accord with the definition on page 30 of \cite{Li}, the
aforementioned Markov process would be a measure-valued branching
process. The uniqueness stated in Theorem \ref{2tm} can be used to
prove that such a process is unique up to modifications. Another
observation is that, in our model, branching is the only
evolutionary act, whereas papers on branching in finite particle
systems, e.g., \cite{BL,BL1,DGL,DKS}, assume more such acts, e.g.,
diffusion in $X$. These and similar generalizations can also be done
in our setting.

\section*{Acknowledgements}
The research of the first named author was financially supported by
National Science Centre, Poland, grant 2017/25/B/ST1/00051, that is
acknowledged by him. The authors are cordially grateful to both
referees for their careful reading of the manuscript, constructive
criticism and valuable and favorable suggestions that helped to
improve the quality of this work.


\begin{thebibliography}{9}



\bibitem{Arendt} W. Arendt, Ch. J. K. Batty, M. Hieber, F.
Neubrander, Vector-valued  Laplace Transforms and Cauchy Problems,
Second Edition, Monographs in Mathematics, Vol. 96, Birkh\"aser,
Basel, 2011.


\bibitem{BL} L. Beznea, O. Lupa\c{s}cu, Measure-valued discrete branching Markov
processes, Trans. Amer. Math. Soc. 368 (2016), 5153--5176.

\bibitem{BL1} L. Beznea, O. Lupa\c{s}cu-Stamate, C. I. Vrabie, Stochastic solutions to evolution equations of non-local branching
processes, Nonlinear Anal. 200 (2020) 112021.

\bibitem{BeR} L. Beznea, M. R\"ockner, From resolvents to cadlag processes through compact
excessive functions and applications to singular SDE on Hilbert
spaces, Bull. Sci. math. 135 (2011) 844--870.

\bibitem{Mich} V. I. Bogachev, N. V. Krylov, M. R\"ockner, S. V. Shaposhnikov,
Fokker-Planck-Kolmogorov Equations. Mathematical Surveys and
Monographs, 207. American Mathematical Society, Providence, RI,
2015.

\bibitem{DV2} D. J. Daley, D. Vere-Jones, An Introduction to the Theory of Point Processes. Vol. II. General Theory and Structure.
Second edition. Probability and its Applications (New York).
Springer-Verlag, New York, 2008.



\bibitem{Cohn} D. L. Cohn, Measure Theory. Second edition. Birkh\"auser Advanced Texts: Basler Lehrb\"ucher.  Birkh\"auser/Springer,
 New York, 2013.


\bibitem{Dawson} D. A. Dawson, Measure-Valued Markov Processes. {\'E}cole d'{\'E}t{\'e} de
Probabilit{\'e}s de Saint-Flour XXI--1991, 1--260, Lecture Notes in
Math., 1541, Springer, Berlin, 1993.


\bibitem{DGL} D. A. Dawson, L. G. Gorostiza, Z. Li, Nonlocal
branching superprocesses and some related models, Acta Applicandae
Mathematicae 74  (2002), 93--112.

\bibitem{DSS}  R. L. Dobrushin, Y. G. Sinai, Y. M. Sukhov, Dynamical systems of statistical
mechanics, in Dynamical Systems II. Encyclopaedia of Mathematical
Sciences, vol 2., Y. G. Sinai, eds, Springer, Berlin, Heidelberg,
1989.


\bibitem{DKS} E. B. Dynkin, S. E. Kuznetsov, A. V. Skorokhod,
Branching measure-valued processes, Probab. Theory Relat, Fields 99
(1994), 55--96.

\bibitem{EK} S. N. Ethier, T. G. Kurtz, Markov Processes: Characterization and Convergence, {\it Wiley, New York,} 1986.





\bibitem{Konar} V. Konarovskyi,  A system of coalescing heavy diffusion particles on the real
line, Ann. Probab. 45 (2017), 3293-3335.

\bibitem{Konar1} V. Konarovskyi, M. von Renesse, Modified massive Arratia flow and Wasserstein
diffusion, Commun. Pure and Appl. Math.  72 (2019), 764--800.

\bibitem{Koz1} Y. Kozitsky, Stochastic branching at the edge: Individual-based modeling of tumor cell
proliferation, J. Evol. Equ.  21 (2021), 2081--2104.

\bibitem{KP} Yu. Kozitsky, K. Pilorz, Random jumps and coalescence in
the continuum: evolution of ststes of an infinite particle system,
Discrete Cont. Dyn-A 40 (2020), 725--752.


\bibitem{Koz2} Y. Kozitsky, M. R\"ockner, A Markov process for an
infinite interacting particle system in the continuum,  Electron. J.
Probab. 26 (2021), article no. 72, 1--53. .


\bibitem{Lenard} A. Lenard, Correletion functions and the uniqueness
of the state in classical statistical mechanics, {\it Comm. Math.
Phys.} {\bf 30} (1973), 35--44.

\bibitem{Li} Z. Li, Measure-Valued Branching Markov Processes, Probability and its
Applications, Springer, Heidelberg Dordrecht London New York, 2011.

\bibitem{MaR} Z.-M. Ma, M. R\"ockner, Introduction to the Theory of
(Nonsymmetric) Dirichlet Forms, Universitext, Springer-Verlag Berlin
Heidelberg, 1992.


\bibitem{Part} K. R. Parthasarathy, Probability Measures on Metric
Spaces,  Probability and Mathematical Statistics, No. 3 Academic
Press, Inc., New York-London 1967.




\bibitem{Simon}  B. Simon, The Statistical Mechanics of Lattice
Gases. I, Princeton University Press, Pronceton, NJ, 1993.


\bibitem{Zessin} H. Zessin, The method of moments for random
measures, Z. Wahrscheinlichkeitstheorie verw. Gebiete. 62 (1983),
395--409.

\end{thebibliography}
\end{document}